\documentclass[twoside,12pt,leqno]{article}
\advance\topmargin by -\headheight\advance\topmargin by -\headsep

\usepackage{amsmath}
\usepackage{amsthm}
\usepackage{amssymb}
\usepackage{amscd}
\usepackage{graphicx}
\usepackage{epsfig}

\numberwithin{equation}{section}
\allowdisplaybreaks

\newtheorem{theorem}{Theorem}[section]
\newtheorem{lemma}{Lemma}[section]

\theoremstyle{definition}
\newtheorem{definition}{Definition}[section]
\newtheorem{example}{Example} [section]

\newtheorem{remark}{Remark}[section]

\begin{document}

\newcommand{\lra}{\longrightarrow}
\newcommand{\C}{\mathcal{C}}
\newcommand{\sM}{\underline{M}}
\newcommand{\X}{\mathcal{X}}
\newcommand{\cY}{\mathcal{Y}}
\newcommand{\Y}{\mathcal{Y}}
\newcommand{\G}{\mathcal{G}}
\newcommand{\cG}{\mathcal{G}}
\newcommand{\cX}{\mathcal{X}}
\newcommand{\sH}{\mathcal{H}}
\newcommand{\cH}{\mathcal{H}}
\newcommand{\bs}{\bar{\mathbf{s}}}
\newcommand{\bt}{\bar{\mathbf{t}}}
\newcommand{\bbbs}{\mathbf{s}}
\newcommand{\bbbt}{\mathbf{t}}
\newcommand{\bm}{\bar{m}}
\newcommand{\be}{\bar{e}}
\newcommand{\bi}{\bar{i}}
\newcommand{\F}{\mathcal{F}}
\newcommand{\R}{\ensuremath{\mathbb R}}

\newcommand {\comment}[1]{{\marginpar{*}\scriptsize{\bf Comments:}\scriptsize{\ #1 \ }}}


\setcounter{page}{1}
\thispagestyle{empty}



\markboth{Hsian-Hua Tseng and Chenchang Zhu}{Integrating Poisson manifolds via stacks}

\label{firstpage}
$ $
\bigskip

\bigskip

\centerline{{\Large  Integrating Poisson manifolds via stacks}}

\bigskip
\bigskip
\centerline{{\large by  Hsian-Hua Tseng and Chenchang Zhu}}

\vspace*{.7cm}

\begin{abstract}
A symplectic groupoid $G.:=(G_1 \rightrightarrows G_0)$ determines a Poisson structure on $G_0$. In this case, we call $G.$ a symplectic groupoid of the Poisson manifold $G_0$. However, not every Poisson manifold $M$ has such a symplectic groupoid. This keeps us away from some desirable goals: for example,  establishing Morita equivalence in the category of all Poisson manifolds. In this paper, we construct symplectic Weinstein groupoids which provide a solution to the above problem (Theorem \ref{main}). More precisely, we show that a symplectic Weinstein groupoid induces a Poisson structure on its base manifold, and that to every Poisson manifold there is an associated symplectic Weinstein groupoid.
\end{abstract}

\pagestyle{myheadings}

\section{Introduction}
The notion of a symplectic groupoid (see \cite{karasev}, \cite{w-poisson}) was introduced in Weinstein's program of quantization of Poisson manifolds. There is an almost 1-1 correspondence between symplectic groupoids and {\em integrable} (to be explained below) Poisson manifolds. This correspondence is closely related to the integration problem of  Lie algebroids, which we now explain.  

Recall that a Lie algebroid over a manifold $M$ is a vector bundle $\pi: A\to M$ with a real Lie bracket structure $[\,,\,]$ on its space of sections $H^0(M,A)$ and a bundle map $\rho :A\to TM$ such that the Leibniz rule $$[X,fY](x)=f(x)[X,Y](x)+(\rho(X)f)(x)Y(x)$$ holds for all $X,Y\in H^0(M,A)$, $f\in C^\infty(M)$ and $x\in M$. The map $\rho$ is called the anchor. It induces a map between $H^0(M,A)$ and the space of global vector fields on $M$, which is a Lie algebra homomorphism.

When $M$ is a point, a Lie algebroid becomes a Lie algebra. A Lie algebra encodes the infinitesimal information of a Lie group. One obtains a Lie algebra from a Lie group by differentiation. One may think the process of obtaining a Lie group from a Lie algebra as a kind of ``integration''. In analogy, a Lie algebroid can be thought of as an infinitesimal version of a Lie groupoid. One can obtain a Lie algebroid from a Lie groupoid by taking invariant vector fields and restricting them to the identity section. The integrability problem of Lie algebroids asks for a reverse process, namely one that associates to a Lie algebroid $A$ a Lie groupoid whose Lie algerboid is $A$. This problem, first formulated  in \cite{pradines} has attracted a lot of attention over time. A solution using local groupoids was also  given in \cite{pradines}, but the global integrating object, which is important for establishing Morita equivalence of Poisson manifolds, was still missing. An important approach to finding such a global object is the use of path spaces. This idea is not new, see \cite{W3} for a nice discussion. We pay particular attention to the recent work \cite{CrF2} of M. Crainic and R. Fernandes and \cite{cafe} of A. Cattaneo and G. Felder. For a Lie algebroid $A$, they study the space of $A$-paths. They are able to give an answer to the integrability problem negatively---not every Lie algebroid can be integrated into a Lie groupoid. From the space of $A$-paths they construct a {\em topological} groupoid and determine equivalent conditions for this groupoid to be a Lie groupoid that integrates the given Lie algebroid $A$. So their work shows that every Lie algebroid can be integrated into a topological groupoid, but in general this topological groupoid doesn't have enough information to recover the Lie algebroid we start with. As conjectured by  Weinstein, one hopes that there are additional structures on this topological groupoid which allow us to recover the Lie algebroid. In \cite{TZ}, the authors find such structures. The key idea is to enlarge the category one works in to the category of {differentiable stacks}. We introduce the notion of {\em Weinstein groupoid} which formalizes the additional structures to put on this topological groupoid. By allowing Weinstein groupoids, we answer the integrability problem positively---every Lie algebroid can be integrated into a Weinstein groupoid.

For a Poisson manifold $M$, there is an associated Lie algebroid $T^*M \to M$. The anchor map $T^*M\to TM$ is given by the contraction with the Poisson bivector, and the Lie bracket is induced by the Poisson bracket $\{, \}$ of $M$, 
\[ [df , dg]:= d\{f, g\} .\]

When $T^*M\to M$ is an integrable Lie algebroid, the Lie groupoid $G_1\rightrightarrows M$ associate to it natually has a multiplicative 2-form $\omega$ on $G_1$. The identity
\[ m^* \omega = pr_1^* \omega + pr_2^* \omega , \]
holds on the composable pairs $G_1 \times_M G_1$, where $m$ is the multiplication and $pr_j$ are projections onto the $j$-th components. It turns out \cite{mw} that the source map of $G_1$ is a Poisson map, that is, the symplectic structure of $G_1$ and the source map recover the Poisson structure  on $M$. In this case, $G_1$ is called a {\em symplectic groupoid} of $M$. Moreover, there is a unquie source-simply connected symplectic groupoid of $M$.

Conversely, to each Poisson manifold $M$, we want to find a symplectic groupoid over it which induces the Poisson structure of $M$. As explained above, one may take the Lie groupoid integrating Lie algebroid $T^*M\to M$. However, even this special kind of Lie algebroids is not always integrable (see for example \cite{bw2} and the references therein). However, if we allow Weinstein groupoids, then there is always such a reverse procedure. In this paper, we construct {\em symplectic Weinstein groupoids} (see Definition \ref{swg}) for every Poisson manifold and prove the analogue of the classical statement above. More precisely, 

\begin{theorem}\label{main} For any symplectic Weinstein groupoid $\cG\rightrightarrows M$, the base manifold $M$ has a unique Poisson structure such that the source map $\bs$ is Poisson. In this case, we call $\cG$ a  symplectic Weinstein groupoid of the Poisson manifold $M$.

 On the other hand, for any Poisson manifold $M$, there are two symplectic groupoid $\cG(T^*M)$ and $\cH(T^*M)$ of $M$.   
\end{theorem}

We also relate the symplectic Weinstein groupoid to the classical symplectic groupoid in the case when $T^*M$ is integrable. 

\begin{theorem}\label{main2} A Poisson manifold $M$ is integrable, i.e. $M$ has an associated symplectic groupoid, if and only if $\cH(T^*M)$ is representable. In this case, $\cH(T^*M)$ the source-simply connected symplectic groupoid integrating $M$.
\end{theorem}

The paper is organized as follows. In Section \ref{diffstack} we recall the notion and properties of differentiable stacks. In Section \ref{weingpd} we discuss the notion of Weinstein groupoids and results in \cite{TZ}. In Section \ref{symplweingpd} we introduce the notion of symplectic Weinstein groupoids and establish a correspondence between symplectic Weinstein groupoids and Poisson manifolds.

{\em Acknowledgements:} We thank K. Behrend, H. Bursztyn,  M. Crainic, T. Graber, D. Metzler, I. Moerdijk, J. Mr{\v{c}}un, D. Salamon, A. Weinstein, P. Xu and M. Zambon for very helpful discussions and suggestions.

\section{Differentiable stacks}\label{diffstack}
In this section we briefly discuss the notion of differentiable stacks. In the past few decades stacks over the category of schemes had be extensively studied in algebraic geometry, especially in connection with moduli problems (see for instance \cite{DM}, \cite{V}, \cite{LMB}, \cite{BEFFGK}). As known in the early days, stacks can also be defined over other categories such as category of topological spaces or smooth manifolds (see for instance \cite{SGA4}, \cite{P}, \cite{BX},\cite{Me}). In this paper we focus on stacks over the category of smooth manifolds, which are called {\em differentiable stacks}. The readers are refered to \cite{P}, \cite{BX} and \cite{Me} for more detailed discussions about differentiable stacks.
\subsection{Definitions}
Let $\C$ be the category of smooth manifolds. A stack over $\C$ is a category fibered in groupiods satisfying two conditions: ``isomorphism is a sheaf'' and ``descend datum is effective''. Both conditions are rather complicated to describe. A precise definition can be found in \cite{BX}, \cite{Me}. See also \cite{F} for an illuminating discussion.

\begin{example}
\hfill
\begin{enumerate}
\item
A manifold $M$ can be viewed as a stack over $\C$. Let $\sM$ denote the following category: the objects are pairs $(S,u)$ where $S$ is a manifold and $u:S\to M$ is a smooth map. A morphism $(S,u)\to (T,v)$ of objects is a smooth map $f:S\to T$ such that $u=v\circ f$. The category $\sM$ is a stack. It contains all the information about the manifold $M$. In this way the notion of stacks generalizes that of manifolds and we identify manifolds with their associated stacks. A stack over $\C$ is called {\em representable} if it is of the form $\sM$ for some manifold $M$.

\item Let $G$ be a Lie group. Recall that the category $BG$ of principal $G$ bundles is defined as follows: the objects are principal bundles $P\to M$ over manifolds. A morphism between two objects $P\to M$ and $P'\to M'$ is a smooth map $M\to M'$ and a $G$-equivariant map $P\to P'$ that covers $M\to M'$. In fact $BG$ is a stack---the classfying stack of $G$-bundles.
\end{enumerate}

\end{example}
Morphisms between stacks are functors between their underlying categories. 
A morphism $f:\X\to \Y$ is a representable submersion if for any morphism $M\to \Y$ from a manifold $M$, the fiber product $\X\times_\Y M$ is representable and the induced morphism $\X\times_\Y M\to M$ is a submersion (between manifolds). If in addition $\X\times_\Y M\to M$ is surjective, then $f$ is called a representable surjective submersion, see \cite{BX}.

\begin{definition}
A differentiable stack $\X$ is a stack over $\C$ together with a representable surjective submersion $\pi:X \to \X$ from a smooth manifold $X$. The morphism $\pi: X\to\X$ is called an {\em atlas} of $\X$. We often abuse notation and call $X$ an atlas of $\X$. Needless to say, atlases are not unique.
\end{definition}
Properties of morphisims between differentiable stacks can be defined by considering pullbacks to atlases. In this way one can define what it means for a morphism to be smooth, \'etale\footnote{In the smooth category, being \'etale means being locally diffeomorphic.}, an immersion, etc. A stack $\X$ is said to be {\em \'etale} if it has an atlas $\pi:X\to \X$ where $\pi$ is \'etale. 
\subsection{Stacks and Groupoids}
Roughly speaking, there is a one-to-one correspondence between differentiable stacks and Lie groupoids. For a differentiable stack $\X$ with an atlas $X_0\to\X$, we obtian a Lie groupoid $X_1:=X_0\times_\X X_0 \rightrightarrows X_0$ where the two maps are projections. This groupoid is called a {groupoid presentation} of $\X$. An \'etale differentiable stack has an \'etale groupoid presentation. Different atlases give different groupoids. But different groupoid presentations of the same stack are {\em Morita equivalent} (see \cite{P}, \cite{BX}, \cite{Me}). Given a groupoid, one can construct a stack (see \cite{V}, \cite{BX}). This process is complicated for a general groupoid. We describe only a special case.
\begin{example}
Consider a Lie group $G$ acting on a manifold $M$. This corresponds to a groupoid $G\times M\rightrightarrows M$ where the two maps are the action and the projection to the second factor. Define a category $[M/G]$ as follows: an object is a principal $G$-bundle $P\to B$ over a manifold $B$ with a $G$-equivariant map $P\to M$. A morphism between two objects $B\leftarrow P\to M$ and $B'\leftarrow P'\to M$ is a pair of a map $B\to B'$ and a $G$-equivarant map $P\to P'$ making all natural diagrams commute. The category $[M/G]$ is in fact a differentable stack with an atlas $M\to [M/G]$.
\end{example}
Differentiable stacks correspond to Morita equivalence classes of Lie groupoids. 1-morphisms between differentiable stacks correspond to what are called {\em Hilsum-Skandalis} morphisms (HS morphisms), see \cite{Mr} and \cite{TZ} for more details.

\section{Weinstein groupoids}\label{weingpd}
\subsection{The Definition}
\begin{definition}[\cite{TZ}]
A Weinstein groupoid over a manifold $M$ consists of the following data: a differentiable stack $\G$, two surjective submersions $\bs,\bt: \G\to M$ (source and target), a map $\bm: \G\times_{\bs,\bt}\G\to \G$ (multiplication), an injective immersion $\be:M\to \G$ (identity section), and an isomorphism $\bi:\G\to\G$ (inverse). These maps are required to satisfy identities\footnote{Those identities are required to hold only up to 2-morphisms.} analogous to those of a groupoid.  
\end{definition}
Roughly speaking, a Weinstein groupoid is a groupoid in the category of differentiable stacks. Let $G$ be the {orbit space} of the stack $\G$, which is a topological space. The data of a Weinstein groupoid induce a {topological} groupoid $G\rightrightarrows M$.
\subsection{The Path Spaces}
We now explain the use of path spaces in the integrability problem. 
\begin{definition}
\hfill
\begin{enumerate}
\item (\cite{CrF2})
Let $\pi:A\to M$ be a Lie algebroid with anchor $\rho: A\to TM$. A $C^1$-map $a: I=[0,1]\to A$ is an $A$-path if the equation $$\rho(a(t))=\frac{d}{dt}(\pi\circ a(t))$$ holds.
\item (\cite{TZ}) Such a map $a: T\to A$ is called an $A_0$-path if, in addition, both $a$ and $\frac{da}{dt}$ vanish on the boundary.
\end{enumerate}
\end{definition}
Denote by $PA$ and $P_0A$ the spaces of $A$- and $A_0$-paths respectively. It's known that $PA$ is a Banach manifold (of infinite dimension) and $P_0A$ is a Banach submanifold of $PA$ (see \cite{CrF2}, Section 1 and \cite{TZ}, Section 2). We next consider the notion of homotopy.
\begin{definition}[\cite{CrF2} \cite{TZ}]
Let $a(\epsilon, t)$ be a family of $A_0$-paths which is $C^2$ in $\epsilon$. Assume that the base paths $\gamma(\epsilon, t):=\rho\circ a(\epsilon,t)$ have fixed end points. For a connection $\nabla$ on $A$, consider the equation 
\begin{equation}\label{hom}
\partial_t b-\partial_\epsilon a=T_\nabla(a,b),\quad b(\epsilon, 0)=0.
\end{equation} Here $T_\nabla$ is the torsion of the connection defined by $T_\nabla(\alpha, \beta)=\nabla_{\rho(\beta)}\alpha-\nabla_{\rho(\alpha)}\beta +[\alpha,\beta].$ Two paths $a_0=a(0,\cdot)$ and $a_1=a(1,\cdot)$ are homotopic if the solution $b(\epsilon,t)$ satisfies $b(\epsilon, 1)=0$.
\end{definition}

\begin{remark}
\hfill
\begin{enumerate}
\item A solution $b(\epsilon,t)$ to (\ref{hom}) doesn't depend on $\nabla$. Therefore the definition makes sense. Furthermore $b(\cdot,t)$ is an $A$-path for each fixed $t$.
\item This definition is analogous to the definition of homotopy of $A$-paths in \cite{CrF2}.
\end{enumerate}
\end{remark}
Homotopies of paths generate foliations $\F$ and $\F_0$ on $PA$ and $P_0A$ respectively. The foliation $\F$ restricts to $\F_0$ on $P_0A$. Now the idea is to consider the {\em monodromy groupoid} (see \cite{MM}) of this foliation: the objects are points in the manifold, and arrows are paths within a leaf (up to homotopies) with fixed end points inside the leaf. Let $Mon(P_0A)\rightrightarrows P_0A$ denote the groupoid associated to the foliation $\F_0$. This groupoid encodes the equivalence relation (i.e. homotopy) of $A_0$-paths. One could think of $Mon(P_0A)$ as the space of homotopies of $A_0$-paths. The two maps from $Mon(P_0A)$ to $P_0A$ assign to each homotopy the two paths at the ends. There are also two maps $P_0A\rightrightarrows M$ which assign to each $A_0$-path its two end points respectively.

Strictly speaking $Mon(P_0A)\rightrightarrows P_0A$ is not a Lie groupoid since both spaces are infinite dimensional. But it is a smooth groupoid in the category of Banach manifolds. Sometimes, to avoid dealing with infinite dimensional issues, we consider a variant $\Gamma\rightrightarrows P$ of this groupoid obtained as follows: $P$ is the disjoint union of an open cover of $P_0A$, and $\Gamma$ is the disjoint union of slices over this cover that are transversal to the foliation $\F_0$, see for instance \cite{TZ} for more details. $\Gamma\rightrightarrows P$ is a finite dimensional Lie groupoid. What's even better is that it's an {\'etale} groupoid (i.e. the source and target are \'etale). The two groupoids $Mon(P_0A)\rightrightarrows P_0A$ and $\Gamma\rightrightarrows P$ are in fact Morita equivalent. Also, there are still two maps $P\rightrightarrows M$.

The next step is clear: We want to consider homotopy equivalence classes of paths and declare that points joined by a homotopy class of paths are equivalent. For this we need to take the ``quotient'' $P_0A/Mon(P_0A)$ and construct a ``groupoid'' $P_0A/Mon(P_0A)\rightrightarrows M$ where the two maps are end-point maps. There are at least two ways to do this. We can take the quotient as a topological space (the topological quotient). Then we obtain a {\em topological} groupoid $P_0A/Mon(P_0A)\rightrightarrows M$ which might not carry any further structure. There is information lost in this process, essentially because the topological quotient remembers only orbits of the equivalence relation $Mon(P_0A)\to P_0A\times P_0A$ given by the groupoid $Mon(P_0A)\rightrightarrows P_0A$ but forgets the finer structures of an orbit.

We can also take the quotient as {stacks}, namely consider the stack associated to the groupoid $Mon(P_0A)\rightrightarrows P_0A$. Given the correspondence between groupoids and stacks, we expect not to lose any information doing this. Denote the stack quotient by $\G:=[P_0A/Mon(P_0A)]$. Since $Mon(P_0A)\rightrightarrows P_0A$ and $\Gamma\rightrightarrows P$ are Morita equivalent, the quotient $[P/\Gamma]$ also equals to $\G$. Since $\Gamma\rightrightarrows P$ is \'etale, $\G$ is an \'etale stack. Moreover, the two maps to $M$ descend to the quotient, giving two maps $\bs, \bt:\G\to M$. There are other maps: By concatenation of paths, we can define a ``multiplication'' $\bm:\G\times_{\bs,\bt}\G\to \G$; by reversing the orientation of a path, we can define an ``inverse'' $\bi:\G\to\G$; by considering constant paths, we can define an ``identity section'' $\be: M\to \G$. These maps are defined in detail in \cite{TZ}. There,  we prove that this makes $\G\rightrightarrows M$ into a Weinstein groupoid.
\begin{theorem}[\cite{TZ}] \label{fromtz}
\hfill
\begin{enumerate}
\item (Lie's third theorem) To each Weinstein groupoid one can associate a Lie algebroid. For every Lie algebroid $A$ there are two natural Weinstein groupoids $\G(A)$ and $\sH(A)$ with Lie algebroid $A$.

\item A Lie algebroid $A$ is integrable in the classical sense if and only if $\sH(A)$ is representable, namely it's a Lie groupoid in the category of manifolds. In this case $\sH(A)$ is the source-simply connected Lie groupoid\footnote{It's called the Weinstein groupoid of A in \cite{CrF2}.} of $A$. 

\item The orbit spaces of $\G(A)$ and $\sH(A)$ (which are topological spaces) are both isomorphic to the universal topological groupoid of $A$ constructed in \cite{CrF2}.

\item Given a Weinstein groupoid $\G$, there is a local groupoid\footnote{It is unique up to isomorphisms near the identity section.} $G_{loc}$ whose Lie algebroid is the same as that of $\G$. 
\end{enumerate}
\end{theorem}

\section{Symplectic Weinstein groupoids}\label{symplweingpd}
In this section we consider the integration problem of Poisson manifolds, namely, the integrability of the Lie algebroid $T^*M\to M$ associated to a Poisson manifold $M$. We introduce the notion of symplectic and Poisson structures on a differentiable stack and apply Theorem \ref{fromtz} to establish a correspondence between Poisson manifolds and what we call symplectic Weinstein groupoids (see Definition \ref{swg}).
\subsection{Symplectic and Poisson Structures}
\begin{definition}
Let $\X$ be a stack over $\C$. The {\em sheaf of differential
$k$-forms} of $\X$  is a contravariant functor $\F^k$ from $\X$ to the category of vector spaces. For every $x\in \X$ over $U\in \C$, define $\F^k(x):=\Omega^k(U)$. For every arrow $y\to x$ over $f:V\to U$, there is a map $\F^k(f):\F^k(x)\to \F^k(y)$ defined by the pullback $f^*: \Omega^k(U)\to\Omega^k(V)$.
\end{definition}
The functor $\F^k$ is in fact a sheaf over $\X$, see \cite{BX} for the definition of sheaves over stacks and the proof of this fact. A {\em differential
$k$-form} $\omega$ on $\X$ is a map that associates to an element
$x\in\X$ over $U$ a section $\omega(x) \in \Omega^k(U)$ such that the
following compatibility condition holds: if there is an arrow $y\to x$
over $f: V\to U$, then $\omega(y)$ is the pull back of $\omega(x)$ by
$f$. Notice that according to this definition, the 0-forms
on $\X$ are simply the maps from $\X$ to $\R$ (viewed as a stack).

There is a simpler interperation when the stack is \'etale:
\begin{lemma}[\cite{z}] \label{etale-form}
Let $\X$ be an \'etale differentiable stack and $G$ an \'etale groupoid
presentation of $\X$. Then there is a 1-1 correspondence between $k$-forms on $\X$ and $G$-invariant $k$-forms on $G_0$.
\end{lemma}
\begin{proof}
A $G$-invariant $k$-form $\omega$ on $G_0$ defines a differential
form on $\cX$ as follows: Given a right $G$-principal
bundle $\pi: P\to U$ with moment map $J: P\to G_0$, the pull back
form $J^*\omega$ is $G$-invariant on $P$. Therefore it induces a
$k$-form $\pi_* J^* \omega$ on $U$ and this is what $P$ associates
to via $\omega$. Notice that we use the fact that $\pi$ is \'etale
to show that a $G$-invariant form is a basic form. On the other
hand, given any $k$-form $\omega$ on $\cX$, consider $\bbbt :G_1\to
G_0$ as a right $G$-principal bundle with moment map $\bbbs: G_1\to G_0$. Then $\omega(G_1)$ is a
$k$-form on $G_0$. Notice that the left multiplication by a certain bisection $ g\cdot : G_1\to G_1$ is a
morphism of $G$-principal bundles. The compatibility
condition of $\omega$ implies that $\omega(G_1)$ is
$G$-invariant.
\end{proof}
\begin{remark} In fact (multi-) vector fields on an \'etale differentiable stacks can also be interpreted as invariant global (multi-) vector fields on an atlas. 
\end{remark} 

\begin{definition}[pull-backs of forms on stacks]
Let $\phi:\Y\to\X$ be a map between stacks and $\omega$ a form
on $\X$. Then $\phi^*\omega$ is a form on $\Y$ defined by
associating to $y\in \Y$ the section $\omega(\phi(y))$.
\end{definition}
\begin{remark} \label{rk-pullback}
We omit here the proof that the above definition is well defined (see for example \cite{z}). Using Lemma \ref{etale-form}, the pull-backs of forms on \'etale differentiable stacks correspond to the
ordinary pull-backs on their \'etale atlases (also see Lemma \ref{poi} for the proof in the case that $\phi$ is $id$).
\end{remark}

By Lemma \ref{etale-form}, we can make the following definition:
\begin{definition} \label{sym-form}
A {\em symplectic form (resp.  Poisson bivector) on an \'etale differentiable stack} $\X$ is a $G$-invariant symplectic  form (resp. Poisson bivector) on
$G_0$, where $G$ is an \'etale presentation of $\X$.
\end{definition}
\begin{remark} Since the source and target maps $\bbbt$ and $\bbbs$ of $G.$ are \'etale, a $G$-invariant form $\omega$ is simply a form satisfying $\bbbs^*\omega=\bbbt^*\omega $.  From Remark \ref{rk-pullback}, a form is symplectic on an \'etale differentiable stack iff it is symplectic on all \'etale presentation.
\end{remark}

\begin{definition}\label{swg}
A Weinstein groupoid $\cG \rightrightarrows M$ is a {\em symplectic Weinstein groupoid} if there is a symplectic form $\omega$ on $\cG$ satisfying the following {\em multiplicative} condition:
\[ \bar{m}^* \omega = pr_1^* \omega + pr_2^* \omega, \]
on $\cG\times_{\bs, M,\bt}\cG$, where $pr_i$ is the projection onto the $i$-th factor.
\end{definition}
\begin{remark}
When $\cG$ is a Lie groupoid, this definition coincides with the definition of symplectic groupoids.
\end{remark}
\subsection{Integrability}
 We will show that after replacing the symplectic groupoid by the symplectic Weinstein groupoid, the correspondence between Poisson manifolds and symplectic groupoids holds for every Poisson manifold. Our approach uses Poisson bracket rather than Poisson bivector.

\begin{lemma} \label{poi}
Give an \'etale differentiable stack $\cX$ with a symplectic form $\omega$, there is a Poisson bracket $\{, \}$ on the algebra $C^{\infty}(\cX)$ of $0$-forms (i.e. smooth functions) on $\cX$.     
\end{lemma}
\begin{proof} Take an \'etale groupoid presentation $G.=(G_1\rightrightarrows G_0)$ of $\cX$ and identify $\cX$ with $BG.$. Then $C^{\infty}(\cX)$ is the set of $G.$-invariant functions $f_G$'s on $G_0$, so it is naturally an algebra. Moreover, $\omega$ appears as a $G.$-invariant symplectic form $\omega_G$ on $G_0$. Therefore, we can define $\{f, g\}_G$---the appearence of $\{f, g\}$ on the presentation $G.$ as  $ \{ f_G, g_G\}_{\omega_G}$, where $\{, \}_{\omega_G}$ is the Poisson bracket defined by $\omega_G$. 

We have to show that the above definition is independent of choices of the \'etale presentations. The groupoid $G.$ is said to be {\em strongly equivalent} to  $H.$ if there is a groupoid morphism $\phi: G. \to H.$ and the H.S. bibundle $E:=G_0\times_{\phi, H_0, \bbbt}H_1$ associated to $\phi$ is a Morita bibundle. If two groupoids are Morita equivalent, they are both strongly equivalent to a third groupoid (see for example \cite{m-orbi}). Hence it suffices to show that if $G.$ is strongly equivelant to $H.$ via $\phi$ (and $E$), then they define the same Poisson bracket. Let $J_G$ and $J_H$ be the moment maps from $E$ to $G_0$ and $H_0$ respectively. Notice that
\[   (G_1 \times_{\bt\circ \phi, H_0} H_1 \to E) \to (G_1\overset{\bbbt}{\lra}  G_0) \]
 is a morphism of $G.$-principal bundles (the first principal bundle is the pull-back of the second via $J_G$). Therefore, we have $$J_G^* \omega_G = J_G^* (\omega(G_1 \overset{\bbbt}{\lra} G_0))=\omega(G_1 \times_{\bbbt\circ \phi, H_0} H_1 \to E).$$
 If we change the presentation from $G.$ to $H.$, the right $G.$-principal bundle  $G_1 \times_{\bt\circ \phi, H_0} H_1 \to E $ transforms via $E$ to an $H.$-principal bundle $(G_1 \times_{H_0}H_1\times_{G_0}E)/G_1 = G_0 \times_{H_0}H_1\times_{H_0} H_1$ over $ E$. On the other hand, this principal bundle is also the pull-back of $H_1\overset{\bbbt}{\lra}H_0$ via $J_H$. So \[J_H^*\omega_H=\omega(G_0 \times_{H_0}H_1\times_{H_0} H_1 \to E)=\omega(G_1 \times_{\bbbt\circ \phi, H_0} H_1 \to E)=J^*_G\omega_G .\] Notice that $J_H^*=J_G^* \phi^*$ and $J_G$ is submersion. We have $\omega_G=\phi^*\omega_H$. Similarly, for functions, we also have $f_G=\phi^*f_H$. Therefore, we  have 
\[ \{ f_H, g_H\}_{\omega_H} = \{ \phi^* f_G, \phi^* g_G \}_{\phi^* \omega_G}= \phi^* (\{ f_G, g_G\}_{\omega_G}). \]
So the Poisson bracket on $\cX$ is well defined.
\end{proof}

Given two  stacks $\cX$ and $\cY$ whose smooth functions form Poisson algebras, a morphism $\cX \to \cY$ is called {\em Poisson} if the induced map $C^{\infty} (\Y)\to C^{\infty} (\X)$ preserves the Poisson brackets. 
 
 Now we can prove Theorem \ref{main} and \ref{main2}.

{\em Proof of Theorem \ref{main}.}
Given a symplectic Weinstein groupoid  $\cG\rightrightarrows M$, we can associate to it a local symplectic groupoid $G_{loc}\rightrightarrows M$. The method is similar to the proof of Theorem \ref{fromtz} (see \cite{TZ}). We recall the idea: Let $G.$ be an \'etale presentation of the stack $\cG$. We devide $M$ into pieces $M_l$ and embed them into $G_0$. Then the local groupoid $G_{loc}$ is obtained by gluing small open neighborhoods $U_l\subset G_0$ of these embedded pieces $M_l$. Then $U_l$ is a local groupoid over $M_l$. The multiplicativity of $\omega$ on $\cG$ implies that the symplectic form $\omega_G|_{U_l}$ is multiplicative. Since the symplectic form $\omega_G$ on the \'etale atlas $G_0$ is invariant under the $G_1$-action and the gluing morphisms are also induced by the $G_1$-action (see Proposition 5.3 in \cite{TZ}),  the multiplicative symplectic forms on the $U_l$'s also glue together to a multiplicative symplectic form on $G_{loc}$.  Therefore, there is a unique Poisson structure $\{, \}_M$ on $M$ such that the source map $\bbbs_{loc}$ of $G_{loc}$ is Poisson.

Notice that the pull-back $\bbbs_{loc}^*f$ of $f \in C^{\infty}(M)$ is locally the same as $\bbbs^*f$ in $C^{\infty}(G_0)$. Since Poisson bracket is a local operation on functions and the Poisson bracket on $\cG$ is defined via  the Poisson bracket on $G_0$, we conclude that the source map $\bbbs: \cG \to M$ is Poisson.     

For the converse, recall that for any Lie algebroid $A$, we can associate two Weinstein groupoids $\cG(A)$ and $\cH(A)$, as discussed in Theorem \ref{fromtz}. We prove the converse statement for $\cG (T^*M)$. The proof for $\cH (T^*M)$
is similar. Let $\omega_c$ be the canonical symplectic form on
$T^*M$. Then according to \cite{cafe}, $\omega_c$ induces a
symplectic form on the path space $P T^*M$. The restriction to the $A$-path space $P_a T^*M$ of this symplectic form has kernel exactly the tangent space of
the foliation $\F$ and is invariant along the foliation. Consider
the \'etale presentation $\Gamma\rightrightarrows P$ of
$\cG (T^*M)$. $P$ is  the transversal of the foliation $\F$,
hence the restricted form is a $\Gamma$-invariant symplectic form.
This form induces a symplectic form $\omega$ on $\cG (T^*M)$.
The multiplicativity of $\omega$ follows from the additivity of
the integrals after examining the definition of $\omega_c$.

It remains to prove that $\bs: \cG(T^*M)\to M$ is Poisson. As shown in the proof of Theorem 1.2 in \cite{TZ}, the local groupoid associated to $\cG(T^*M)$ is exactly the symplectic local groupoid associated to $M$ in \cite{cf2}. An argument analogous to the above shows that $\bs$ is a Poisson map.   

{\em Proof of Theorem \ref{main2}.} It follows from Theorem \ref{main} and Theorem \ref{fromtz}.


\noindent 
Hsian-Hua Tseng\\
Department of Mathematics\\
University of California at Berkeley\\
Berkeley, CA 94720, USA\\
hhtseng@math.berkeley.edu\\

\noindent
Chenchang Zhu\\
Departement Mathematik\\
ETH Zentrum\\
R\"amistrasse 101, 8092 Z\"urich, Switzerland\\
zhu@math.ethz.ch


\end{document}